\documentclass[fleqn,twoside]{article}
\usepackage{espcrc1}

\newcommand{\ttbs}{\char'134}
\newcommand{\AmS}{{\protect\the\textfont2
  A\kern-.1667em\lower.5ex\hbox{M}\kern-.125emS}}

\title{Note on Combinatorial Engineering Frameworks
 for Hierarchical Modular Systems}

\author{Mark Sh. Levin
\thanks{
 Mark Sh. Levin:~
 Inst. for Inform. Transmission Problems,
 Russian Academy of Sciences;
  http://www.mslevin.iitp.ru;
 email: mslevin@acm.org
  }
  }

\begin{document}

\maketitle

\begin{abstract}
 The paper briefly describes
 a basic set of special combinatorial engineering frameworks
 for solving complex problems in the field of hierarchical modular
 systems.
 The frameworks consist
 of combinatorial problems (and corresponding models),
 which are interconnected/linked (e.g., by preference relation).
 Mainly, hierarchical morphological system model is used.
 The list of basic
  standard combinatorial engineering
 (technological)
  frameworks is the following:
 (1) design of system hierarchical model,
 (2) combinatorial synthesis
 ('bottom-up' process for system design),
 (3) system evaluation,
 (4) detection of system bottlenecks,
 (5) system improvement (re-design, upgrade),
 (6) multi-stage design (design of system trajectory),
 (7) combinatorial modeling of system evolution/development
 and system forecasting.
 The combinatorial engineering frameworks are targeted to maintenance
 of some system life cycle stages.
 The list of main underlaying combinatorial optimization problems
 involves the following:
 knapsack problem, multiple-choice problem, assignment problem,
 spanning trees, morphological clique problem.

~~


\vspace{1pc}
\end{abstract}

\maketitle

\newcounter{cms}
\setlength{\unitlength}{1mm}

\section{Introduction}

 The frame  approach
  for representing  knowledge
 (i.e., collection of frames are linked together into frame-system)
 has been suggested by Marvin Minsky
 (e.g., \cite{minsky75}).
 In this approach, the frame corresponds to a data structure.
 In general,
 it is possible to consider the following three-component system:
 initial data/information,
 problem(s) (and corresponding models),
 and
 algorithm (or interactive procedure).
 For many complex applied problems,
 it is reasonable to examine special composite frameworks
 (i.e., composite solving schemes)
 consisting of problems (and corresponding models),
 which are interconnected/linked (e.g., by preference relation).
 For example,
 a basic simplified framework for data processing
 can be described as follows:

 (a) analysis of input data/information and preliminary processing;

 (b) processing; and

 (c) analysis of results.

 Another example of a framework is well-known in decision
 making.
 Herbert Simon has suggested
 his framework
 for rational decision making (choice problem)
 (e.g., \cite{sim76}):
 (i) the identification and listing of all the alternatives,
 (ii) determination of all the consequences resulting from each of the
 alternatives, and
 (iii) the comparison of the accuracy and efficiency of each of
 theses sets of consequences.
%
%
  A modified version of this decision making framework is the following:

 {\it Stage 1.} Analysis of the examined system/process, extraction of the problem.

 {\it Stage 2.} Problem structuring:

 (2.1.) generation of alternatives,

  (2.2.) generation of criteria and a scale for each criterion.

 {\it Stage 3.} Obtaining the initial information
 (estimates of the alternatives, preferences over the alternatives).

 {\it Stage 4.} Solving process to obtain the decision(s).

 {\it Stage 5.} Analysis of the obtained decision(s).

 On our opinion,
   there exists a trend to design, describe,
   and use a set of basic typical
 engineering (technological) frameworks
 (i.e., typical composite combinatorial solving schemes),
 which can be considered as basic standard blocks
 in systems research/design and
 in systems education
 (engineering, computer science, applied mathematics).

 In recent decades, modular approaches have been used in all engineering
 domains.
 Thus,  many systems can be designed from basic standard modules
 (e.g., software engineering, computer engineering, information
 engineering, method engineering, protocol engineering).
 Evidently, special combinatorial methods have to be
 studied and applied for system analysis/design
 at all stages of the system life cycles.
 The methods can have the following structure types:
 series, parallel, series-parallel,
 cascade-like.
 Here, two basic problems are very important:
 (1) partitioning the initial problem
 (or partitioning the examined system),
 (2) aggregation of solutions.
%

 This note contains our attempt to describe
 the basic set of
 typical combinatorial engineering frameworks
 for hierarchical modular systems
 (with hierarchical structures).
 This approach is based on
  the following five-layer
  architecture
 (a modification of the architecture from \cite{lev11ADES}):

 Layer 1. Basic combinatorial optimization problems
 (e.g., knapsack problem, multiple choice problem, multicriteria ranking/selection,
 clustering, minimum spanning tree, Steiner tree, clique problem).

 Layer 2. Complex (e.g., multicriteria) combinatorial optimization
 problems
 (e.g., multicriteria
 knapsack problem, multicriteria multiple choice problem,
 multicriteria Steiner tree, morphological clique problem,
 design of multi-layer network topology).

 Layer 3. Basic support frameworks
 (e.g., hierarchical design, aggregation of structures,
 restructuring of knapsack problem, restructuring of multiple choice problem).

 Layer 4. Combinatorial
 engineering frameworks (consisting of a set
 of
 linked combinatorial problems)
 (e.g., hierarchical system modeling,
  design,
 evaluation,
 detection of bottlenecks,
  improvement, trajectory design,
 combinatorial evolution and forecasting).

 Layer 5. Applied combinatorial
 engineering frameworks
 (e.g., modeling, design and improvement of  components/parts for
 various applied systems).


 Here,
 the basic set of standard combinatorial engineering
 (technological)
  frameworks for modular systems
 are
  described
 while taking into account the system life cycles
 (i.e., layer 4 above)
 \cite{lev06,levprob07}:

 (1) design of system hierarchical model,

 (2) combinatorial synthesis
 ('bottom-up' process for system design),

 (3) system evaluation,

 (4) detection of system bottlenecks
 (e.g., by system elements, by  compatibility
 of system elements, by system structure),

 (5) system improvement (re-design, upgrade),

 (6) multi-stage design (design of system trajectory),

 (7) combinatorial modeling of system evolution/development
 and system forecasting.

 The above-mentioned combinatorial synthesis
 is used as a basic combinatorial engineering framework
 This framework is based on two approaches:
 (1) our modification of morphological analysis
 (Hierarchical Multicriteria Morphological Design HMMD)
 (e.g., \cite{lev98,lev06,lev09,lev12morph,lev12a})
 or
 (2) multiple choice problem
 (e.g., \cite{keller04,lev09,lev12morph,lev12a}).
 HMMD and multiple choice problem
 are used with two kinds of estimates for design
 alternatives of system components/parts:
 (i) ordinal estimates
  (e.g., \cite{lev98,lev06,lev09,lev12morph}),
 (ii) interval multiset estimates
  (e.g., \cite{lev12a,lev13}).

 The material contains the author viewpoint
 to maintenance of modular hierarchical systems
 (i.e., physical systems, software, organizational systems,
 plans, composite solving strategy,
  system requirements, standards,
 communication protocols).

 Some author's system applications are pointed out, for example:
 electronic shopping, Web-based system,
 decision support system,
 modular software,
 composite strategy for multicriteria ranking,
 integrated security system, telemetry system,
 two-floor building,
 control system for smart homes,
 system of political management,
 concrete technology,
 medical treatment, immunoassay technology,
 wireless sensor, communication protocol, and standard for
 multimedia information.

\section{Preliminaries}

 In recent decades, standardization became to be a real basis
 for extensive examination of modular systems in all
  domains of engineering and information technology
  including all stages of system life cycle
  (e.g., system design, system maintenance, system testing,
  etc.).
 On the other hand,
 hierarchical approaches are
 power tools for modeling,
 analysis, and design of various systems.
 Fig. 1 depicts the considered domain as
 system applications,
 system hierarchical structure,
 and
 the basic set of
 combinatorial
 engineering frameworks
 for modular systems.
 A ``two-dimensional'' domain
 for relation between problem(s)/ model(s)
 and algorithm(s)/solving frameworks
 is shown in Fig. 2.
 This representation illustrates
 two system directions as an extension of traditional
 pair ``problem/model - algorithm/solving procedure''.
%
 Our approach is based
 on typical combinatorial engineering frameworks
 as k-problem/k-model frameworks for modular systems.
 Thus, hierarchical modular system models
 and the above-mentioned combinatorial engineering frameworks
 are a fundamental
 for problem structuring and solving
 in real-world system applications.

\begin{center}
\begin{picture}(113,73)
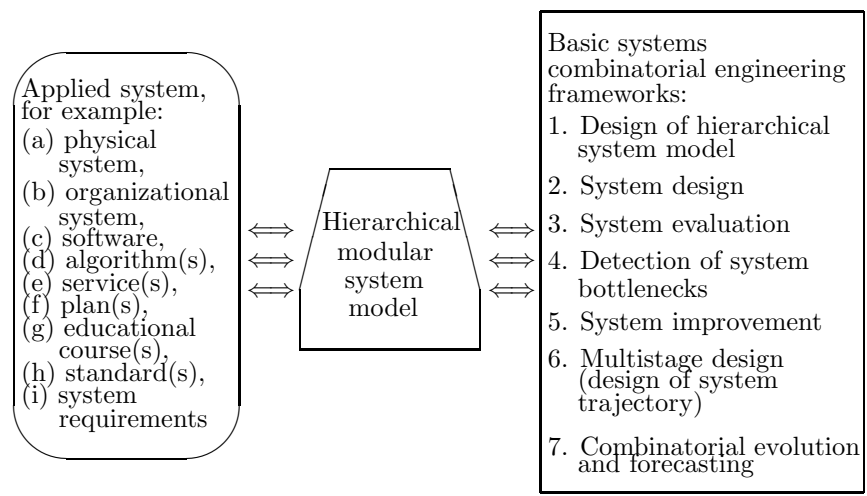
\put(07,00){\makebox(0,0)[bl] {Fig. 1. System - hierarchical model
- engineering frameworks}}


\put(01,60){\makebox(0,0)[bl]{Applied system,}}
\put(01,57){\makebox(0,0)[bl]{for example:}}

\put(01,53){\makebox(0,0)[bl]{(a) physical }}
\put(06,50){\makebox(0,0)[bl]{system,}}

\put(01,46){\makebox(0,0)[bl]{(b) organizational}}
\put(06,43){\makebox(0,0)[bl]{system,}}

\put(01,40){\makebox(0,0)[bl]{(c) software,}}

\put(01,37){\makebox(0,0)[bl]{(d) algorithm(s),}}
\put(01,34){\makebox(0,0)[bl]{(e) service(s),}}

\put(01,31){\makebox(0,0)[bl]{(f) plan(s),}}
\put(01,28){\makebox(0,0)[bl]{(g) educational}}
\put(06,25){\makebox(0,0)[bl]{course(s),}}

\put(01,22){\makebox(0,0)[bl]{(h) standard(s),}}
\put(01,19){\makebox(0,0)[bl]{(i) system }}
\put(06,16){\makebox(0,0)[bl]{requirements}}

\put(15,39.5){\oval(30,54)}


\put(31,42){\makebox(0,0)[bl]{\(\Longleftrightarrow\)}}
\put(31,38){\makebox(0,0)[bl]{\(\Longleftrightarrow\)}}
\put(31,34){\makebox(0,0)[bl]{\(\Longleftrightarrow\)}}


\put(41,43){\makebox(0,0)[bl]{Hierarchical}}
\put(43,39){\makebox(0,0)[bl]{modular}}
\put(44,34.5){\makebox(0,0)[bl]{system}}
\put(44.5,31.5){\makebox(0,0)[bl]{model}}

\put(42,51){\line(1,0){16}}

\put(38,35){\line(1,4){4}} \put(62,35){\line(-1,4){4}}

\put(38,27){\line(0,1){08}} \put(62,27){\line(0,1){08}}

\put(38,27){\line(1,0){24}}

\put(63,42){\makebox(0,0)[bl]{\(\Longleftrightarrow\)}}
\put(63,38){\makebox(0,0)[bl]{\(\Longleftrightarrow\)}}
\put(63,34){\makebox(0,0)[bl]{\(\Longleftrightarrow\)}}


\put(71,66){\makebox(0,0)[bl]{Basic systems }}
\put(71,62.6){\makebox(0,0)[bl]{combinatorial engineering}}
\put(71,60){\makebox(0,0)[bl]{frameworks:}}

\put(71,55){\makebox(0,0)[bl]{1. Design of hierarchical}}
\put(75,52){\makebox(0,0)[bl]{system model}}

\put(71,47){\makebox(0,0)[bl]{2. System design}}

\put(71,42){\makebox(0,0)[bl]{3. System evaluation }}

\put(71,37){\makebox(0,0)[bl]{4. Detection of system}}
\put(75,34){\makebox(0,0)[bl]{bottlenecks}}

\put(71,29){\makebox(0,0)[bl]{5. System improvement}}

\put(71,24){\makebox(0,0)[bl]{6. Multistage design }}
\put(75,21){\makebox(0,0)[bl]{(design of system}}
\put(75,18){\makebox(0,0)[bl]{trajectory)}}

\put(71,13){\makebox(0,0)[bl]{7. Combinatorial evolution}}
\put(75,10){\makebox(0,0)[bl]{and forecasting}}

\put(70,08){\line(1,0){43}} \put(70,72){\line(1,0){43}}
\put(70,08){\line(0,1){64}} \put(113,08){\line(0,1){64}}

\end{picture}
\end{center}


\begin{center}
\begin{picture}(117,83.5)

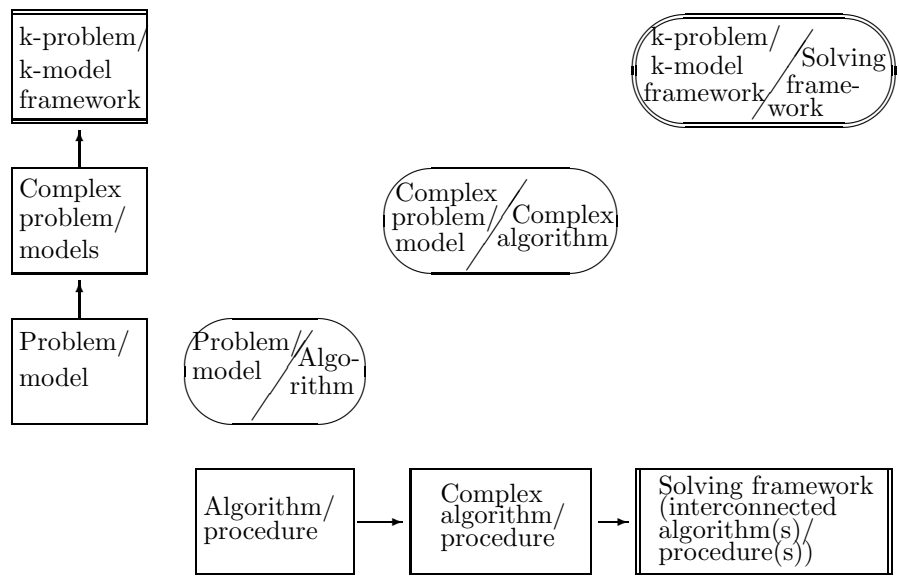
\put(08,00){\makebox(0,0)[bl]{Fig. 2. Domain ``Problem/model -
Algorithm/solving framework''}}


\put(24.5,06){\line(1,0){21}} \put(24.5,20){\line(1,0){21}}
\put(24.5,06){\line(0,1){14}} \put(45.5,06){\line(0,1){14}}

\put(25.5,13){\makebox(0,0)[bl]{Algorithm/}}
\put(25.5,10){\makebox(0,0)[bl]{procedure}}

\put(46,13){\vector(1,0){6}}


\put(53,06){\line(1,0){24}} \put(53,20){\line(1,0){24}}
\put(53,06){\line(0,1){14}} \put(77,06){\line(0,1){14}}

\put(57,15){\makebox(0,0)[bl]{Complex}}
\put(57,12){\makebox(0,0)[bl]{algorithm/}}
\put(57,09){\makebox(0,0)[bl]{procedure}}

\put(78,13){\vector(1,0){4}}


\put(83,06){\line(1,0){34}} \put(83,20){\line(1,0){34}}
\put(83,06){\line(0,1){14}} \put(117,06){\line(0,1){14}}

\put(83.5,06){\line(0,1){14}} \put(116.5,06){\line(0,1){14}}

\put(86,16){\makebox(0,0)[bl]{Solving framework}}
\put(86,13){\makebox(0,0)[bl]{(interconnected}}
\put(86,10){\makebox(0,0)[bl]{algorithm(s)/}}
\put(86,07){\makebox(0,0)[bl]{procedure(s))}}


\put(09,40){\vector(0,1){5}}

\put(00,26){\line(1,0){18}} \put(00,40){\line(1,0){18}}
\put(00,26){\line(0,1){14}} \put(18,26){\line(0,1){14}}

\put(01,35){\makebox(0,0)[bl]{Problem/}}
\put(01,31){\makebox(0,0)[bl]{model}}


\put(35,33){\oval(24,14)}

\put(32,26.5){\line(2,3){08}}

\put(24,35){\makebox(0,0)[bl]{Problem/}}
\put(24,32){\makebox(0,0)[bl]{model}}

\put(38,33){\makebox(0,0)[bl]{Algo-}}
\put(37,30){\makebox(0,0)[bl]{rithm}}


\put(09,60){\vector(0,1){5}}

\put(00,46){\line(1,0){18}} \put(00,60){\line(1,0){18}}
\put(00,46){\line(0,1){14}} \put(18,46){\line(0,1){14}}

\put(01,55.5){\makebox(0,0)[bl]{Complex }}
\put(01,51.5){\makebox(0,0)[bl]{problem/}}
\put(01,48){\makebox(0,0)[bl]{models}}


\put(65,53){\oval(31,14)}

\put(60.5,46.5){\line(2,3){08}}

\put(51,55){\makebox(0,0)[bl]{Complex}}
\put(50.5,51.6){\makebox(0,0)[bl]{problem/}}
\put(51,49){\makebox(0,0)[bl]{model}}

\put(66.5,52){\makebox(0,0)[bl]{Complex}}
\put(64.5,49){\makebox(0,0)[bl]{algorithm}}


\put(00,66.5){\line(1,0){18}} \put(00,80.5){\line(1,0){18}}

\put(00,66){\line(1,0){18}} \put(00,81){\line(1,0){18}}
\put(00,66){\line(0,1){15}} \put(18,66){\line(0,1){15}}

\put(01,75.5){\makebox(0,0)[bl]{k-problem/ }}
\put(01,72){\makebox(0,0)[bl]{k-model}}
\put(01,68){\makebox(0,0)[bl]{framework}}


\put(100,73){\oval(34,14)} \put(100,73){\oval(35,15)}

\put(98.5,66.5){\line(2,3){08}}

\put(85,75.5){\makebox(0,0)[bl]{k-problem/}}
\put(85,72.5){\makebox(0,0)[bl]{k-model}}
\put(84,69){\makebox(0,0)[bl]{framework}}

\put(105,73){\makebox(0,0)[bl]{Solving }}
\put(103,70){\makebox(0,0)[bl]{frame-}}
\put(100.5,67){\makebox(0,0)[bl]{work}}


\end{picture}
\end{center}


\section{Towards Hierarchies}

 Hierarchies play a central role
 in system science,
 in engineering, and
 in computer science
 (e.g., \cite{gar79,haimes90,harel87,kemp08,knuth68,moses10,tanen06}).
 Generally, it is reasonable to point out some
 basic types  of hierarchies
 (e.g., \cite{gar79,haimes90,knuth68,lev12hier}):
 (1) various kinds of trees
 (e.g.,  \cite{gar79,knuth68});
 (2) organic hierarchy
 (i.e., with organic interconnection among children-vertices)
  \cite{conant74};
 (3) ``morphological hierarchy''
  (e.g., \cite{lev98,lev06,lev11agg,lev12morph}); and
 (4) multi-layer structures (Fig. 3)
  (e.g., \cite{bel97,lev11ADES,lev12hier,mes70,obr10,tanen06}).

\begin{center}
\begin{picture}(55,51)

\put(03,00){\makebox(0,0)[bl]{Fig. 3. Multilayer structure}}


\put(39,46){\makebox(0,0)[bl]{Top}}
\put(39,43){\makebox(0,0)[bl]{layer}}

\put(20,46){\oval(32,6)}

\put(09,46){\circle*{1.8}} \put(13,46){\circle*{1.8}}
\put(17.5,45.5){\makebox(0,0)[bl]{{\bf .~.~.}}}

\put(27,46){\circle*{1.8}}  \put(31,46){\circle*{1.8}}


\put(09,46){\vector(0,-1){06}} \put(09,46){\vector(-1,-2){03}}

\put(13,46){\vector(0,-1){06}}

\put(27,46){\vector(0,-1){06}} \put(27,46){\vector(-1,-2){03}}

\put(31,46){\vector(0,-1){06}} \put(31,46){\vector(-1,-2){03}}

\put(31,46){\vector(1,-2){03}}


\put(40,36){\makebox(0,0)[bl]{Intermediate}}
\put(40,33){\makebox(0,0)[bl]{layer}}

\put(20,36){\oval(36,6)}

\put(07,36){\circle*{1.5}} \put(12,36){\circle*{1.5}}
\put(17.5,35.5){\makebox(0,0)[bl]{{\bf .~.~.}}}

\put(28,36){\circle*{1.5}}  \put(33,36){\circle*{1.5}}


\put(07,36){\vector(0,-1){06}} \put(07,36){\vector(-1,-2){03}}
\put(07,36){\vector(1,-2){03}}

\put(12,36){\vector(0,-1){06}} \put(12,36){\vector(1,-2){03}}

\put(28,36){\vector(0,-1){06}} \put(28,36){\vector(-1,-2){03}}

\put(33,36){\vector(0,-1){06}} \put(33,36){\vector(-1,-2){03}}
\put(33,36){\vector(1,-2){03}}



\put(15.8,26.5){\makebox(0,0)[bl]{{\bf . ~. ~.}}}


\put(40,20){\makebox(0,0)[bl]{Intermediate}}
\put(40,17){\makebox(0,0)[bl]{layer}}

\put(20,20){\oval(36,6)}

\put(07,20){\circle*{1.5}} \put(12,20){\circle*{1.5}}
\put(17.5,19.5){\makebox(0,0)[bl]{{\bf .~.~.}}}

\put(28,20){\circle*{1.5}}  \put(33,20){\circle*{1.5}}


\put(07,20){\vector(0,-1){06}} \put(07,20){\vector(-1,-2){03}}
\put(07,20){\vector(1,-2){03}}

\put(12,20){\vector(0,-1){06}} \put(12,20){\vector(1,-2){03}}

\put(28,20){\vector(0,-1){06}} \put(28,20){\vector(-1,-2){03}}

\put(33,20){\vector(0,-1){06}} \put(33,20){\vector(-1,-2){03}}
\put(33,20){\vector(1,-2){03}}



\put(41,10){\makebox(0,0)[bl]{Bottom}}
\put(41,07){\makebox(0,0)[bl]{layer}}

\put(20,10){\oval(40,6)}

\put(05,10){\circle*{1}} \put(10,10){\circle*{1}}
\put(17.5,09.5){\makebox(0,0)[bl]{{\bf .~.~.}}}

\put(30,10){\circle*{1}}  \put(35,10){\circle*{1}}

\end{picture}
\end{center}

 A survey of design methoids for hierarchical multi-layer
 structures is presented in \cite{lev12hier}.
 In the case of trees, spanning tree problems are mainly used to
 design the tree-like hierarchy
 (e.g., minimum spanning tree problems,
 Steiner tree problems, maximum leaf spanning tree problem).
 Some methods for design of 'optimal' organizational hierarchies
 (mainly: trees)
 are examined as well
 (e.g., \cite{bal06,mishin07}).
 On the other hand,
 the following methods are used:
 various expert procedures,
 clustering
 (e.g., hierarchical clustering),
 ontology-based approaches,
 Approaches to design of hierarchical networks are based on
 special combinatorial optimization problems
 (e.g., \cite{bala94,current86,obr10,pirkul91}).

 A general design framework for multi-layer structures
 can be considered as the following
 \cite{lev12hier}:

 {\it Stage 1.} Partitioning the initial set of nodes
 into layer subsets.

 {\it Stage 2.} Design of a topology at each layer.

 {\it Stage 3.} Connection between nodes of neighbor layers.

 On the other hand, it is possible to present basic topology design problems
 from the viewpoint of
 multi-layer topology design (Table 1).

\begin{center}
\begin{picture}(105,100)

\put(12,96){\makebox(0,0)[bl]{Table 1. Design
 approaches for multi-layer network}}

\put(00,00){\line(1,0){105}} \put(00,88){\line(1,0){105}}
\put(00,94){\line(1,0){105}}

\put(00,0){\line(0,1){94}} \put(30,0){\line(0,1){94}}
\put(65,0){\line(0,1){94}} \put(105,0){\line(0,1){94}}

\put(01,89.4){\makebox(0,0)[bl]{Basic problem }}
\put(31,89.4){\makebox(0,0)[bl]{Layers}}
\put(66,89.4){\makebox(0,0)[bl]{Additional problem(s)}}


\put(01,83){\makebox(0,0)[bl]{1.Spanning tree }}
\put(04,79){\makebox(0,0)[bl]{(forest)}}

\put(31,83){\makebox(0,0)[bl]{1.Root(s)}}
\put(31,79){\makebox(0,0)[bl]{2.Transmission nodes}}
\put(31,74){\makebox(0,0)[bl]{3.Leaf nodes}}


\put(01,69){\makebox(0,0)[bl]{2.Steiner tree }}
\put(04,65){\makebox(0,0)[bl]{(forest)}}

\put(31,69){\makebox(0,0)[bl]{1.Root(s)}}
\put(31,65){\makebox(0,0)[bl]{2.Transmission nodes}}
\put(31,61){\makebox(0,0)[bl]{3.Leaf nodes}}

\put(66,69){\makebox(0,0)[bl]{1.Selection/positioning }}
\put(66,65){\makebox(0,0)[bl]{of Steiner nodes}}


\put(01,56){\makebox(0,0)[bl]{3.Maximum leaf }}
\put(04,52){\makebox(0,0)[bl]{nodes}}

\put(31,56){\makebox(0,0)[bl]{1.Root}}
\put(31,52){\makebox(0,0)[bl]{2.Transmission nodes}}
\put(31,48){\makebox(0,0)[bl]{3.Leaf nodes}}

\put(66,56){\makebox(0,0)[bl]{1.Topology over }}
\put(66,52){\makebox(0,0)[bl]{transmission nodes}}


\put(01,43){\makebox(0,0)[bl]{4.Connected  }}
\put(04,39){\makebox(0,0)[bl]{dominating set}}

\put(31,43){\makebox(0,0)[bl]{1.Dominating set}}
\put(31,39){\makebox(0,0)[bl]{2.Leaf nodes}}

\put(66,43){\makebox(0,0)[bl]{1.Topology for }}
\put(66,39){\makebox(0,0)[bl]{dominating set}}
\put(66,35){\makebox(0,0)[bl]{(e.g., tree, path, paths,}}
\put(66,31){\makebox(0,0)[bl]{star, cycle}}


\put(01,26){\makebox(0,0)[bl]{5.Clustering }}

\put(31,26){\makebox(0,0)[bl]{1.Cluster heads}}
\put(31,22){\makebox(0,0)[bl]{2.Cluster nodes}}

\put(66,26){\makebox(0,0)[bl]{1.Selection/assignment}}
\put(66,22.5){\makebox(0,0)[bl]{of cluster heads}}

\put(66,18){\makebox(0,0)[bl]{2.Topology for set}}
\put(66,14.5){\makebox(0,0)[bl]{of cluster heads}}

\put(66,10){\makebox(0,0)[bl]{3.Topology for set}}
\put(66,06.5){\makebox(0,0)[bl]{of cluster nodes}}
\put(66,02){\makebox(0,0)[bl]{(for each cluster)}}

\end{picture}
\end{center}

%
 Fig. 4 depicts an example of multi-layer structure as
 a six-layer communication network.

\begin{center}
\begin{picture}(88,66)
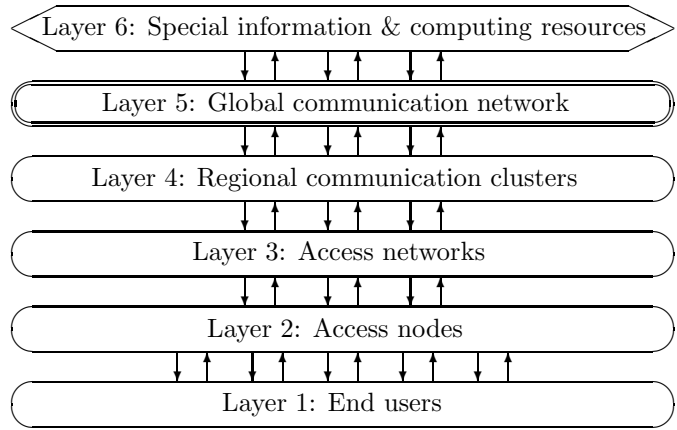
\put(05,0){\makebox(0,0)[bl]{Fig. 4. Six-layer model of
communication network}}



\put(00,59){\line(2,1){06}} \put(00,59){\line(2,-1){06}}

\put(88,59){\line(-2,1){06}} \put(88,59){\line(-2,-1){06}}

\put(06,56){\line(1,0){76}} \put(06,62){\line(1,0){76}}

\put(4,57.5){\makebox(0,0)[bl]{Layer \(6\): Special information \&
 computing resources}}

\put(31,56){\vector(0,-1){4}} \put(35,52){\vector(0,1){4}}
\put(42,56){\vector(0,-1){4}} \put(46,52){\vector(0,1){4}}
\put(53,56){\vector(0,-1){4}} \put(57,52){\vector(0,1){4}}


\put(12,47.3){\makebox(0,0)[bl]{Layer \(5\): Global communication
network}}

\put(44,49){\oval(88,6)} \put(44,49){\oval(87,05)}

\put(31,46){\vector(0,-1){4}} \put(35,42){\vector(0,1){4}}
\put(42,46){\vector(0,-1){4}} \put(46,42){\vector(0,1){4}}
\put(53,46){\vector(0,-1){4}} \put(57,42){\vector(0,1){4}}


\put(010.6,37.4){\makebox(0,0)[bl]{Layer \(4\):  Regional
communication clusters}}

\put(44,39){\oval(88,06)}

\put(31,36){\vector(0,-1){4}} \put(35,32){\vector(0,1){4}}
\put(42,36){\vector(0,-1){4}} \put(46,32){\vector(0,1){4}}
\put(53,36){\vector(0,-1){4}} \put(57,32){\vector(0,1){4}}


\put(24,27.5){\makebox(0,0)[bl]{Layer \(3\): Access networks}}

\put(44,29){\oval(88,06)}

\put(31,26){\vector(0,-1){4}} \put(35,22){\vector(0,1){4}}
\put(42,26){\vector(0,-1){4}} \put(46,22){\vector(0,1){4}}
\put(53,26){\vector(0,-1){4}} \put(57,22){\vector(0,1){4}}


\put(26,17.3){\makebox(0,0)[bl]{Layer \(2\): Access nodes}}

\put(44,19){\oval(88,6)}

\put(22,16){\vector(0,-1){4}} \put(26,12){\vector(0,1){4}}
\put(32,16){\vector(0,-1){4}} \put(36,12){\vector(0,1){4}}
\put(42,16){\vector(0,-1){4}} \put(46,12){\vector(0,1){4}}
\put(52,16){\vector(0,-1){4}} \put(56,12){\vector(0,1){4}}
\put(62,16){\vector(0,-1){4}} \put(66,12){\vector(0,1){4}}


\put(28,07.4){\makebox(0,0)[bl]{Layer \(1\):  End users}}

\put(44,09){\oval(88,06)}


\end{picture}
\end{center}

 In our research projects,
 the above-mentioned special morphological hierarchy
 for system modeling is used (Fig. 5)
 (e.g., \cite{lev98,lev06,lev11agg,lev12morph}.

\begin{center}
\begin{picture}(74,53)
\put(00,00){\makebox(0,0)[bl]{Fig. 5.
 ``Morphological'' system hierarchy
 \cite{lev11agg,lev11ADES}}}


\put(37,51){\circle*{2.5}}


\put(23.6,43){\makebox(0,0)[bl]{System hierarchy}}
\put(21.5,40){\makebox(0,0)[bl]{(tree-like structure)}}


\put(01,34){\makebox(0,0)[bl]{Leaf vertices}}

\put(04,33){\line(0,-1){8}}

\put(05,33){\line(2,-3){5.4}}

\put(6,33){\line(3,-1){23}}


\put(24,36){\makebox(0,0)[bl]{Alternatives for}}
\put(27.5,33.5){\makebox(0,0)[bl]{leaf vertices}}

\put(30,33){\line(-1,-1){14}} \put(43,33){\line(1,-1){15}}


\put(52,36){\makebox(0,0)[bl]{Compatibility}}
\put(52,33){\makebox(0,0)[bl]{among}}
\put(52,30.5){\makebox(0,0)[bl]{alternatives}}

\put(60.4,30){\line(-1,-1){20.5}}


\put(00,23){\line(1,0){74}}

\put(00,39){\line(3,1){37}} \put(74,39){\line(-3,1){37}}

\put(00,23){\line(0,1){16}} \put(74,23){\line(0,1){16}}

\put(02,24.5){\makebox(0,0)[bl]{\(1\)}}

\put(03,23){\circle*{1.8}} \put(03,18){\oval(5,8)}

\put(12,24.5){\makebox(0,0)[bl]{\(2\)}}

\put(13,23){\circle*{1.8}} \put(13,18){\oval(5,8)}

\put(29.5,24){\makebox(0,0)[bl]{\(\tau-1\)}}

\put(32,23){\circle*{1.8}} \put(32,18){\oval(5,8)}

\put(41,24.5){\makebox(0,0)[bl]{\(\tau\)}}

\put(42,23){\circle*{1.8}} \put(42,18){\oval(5,8)}

\put(57,24.5){\makebox(0,0)[bl]{\(m-1\)}}

\put(61,23){\circle*{1.8}} \put(61,18){\oval(5,8)}

\put(70,24.5){\makebox(0,0)[bl]{\(m\)}}

\put(71,23){\circle*{1.8}} \put(71,18){\oval(5,8)}

\put(19.5,17.5){\makebox(0,0)[bl]{. . .}}
\put(48.5,17.5){\makebox(0,0)[bl]{. . .}}

\put(19.5,9){\makebox(0,0)[bl]{. . .}}
\put(48.5,9){\makebox(0,0)[bl]{. . .}}
\put(06,06){\line(1,0){04}} \put(06,13){\line(1,0){04}}
\put(06,06){\line(0,1){7}} \put(10,06){\line(0,1){7}}

\put(07,06){\line(0,1){7}} \put(08,06){\line(0,1){7}}
\put(09,06){\line(0,1){7}}

\put(35,06){\line(1,0){04}} \put(35,13){\line(1,0){04}}
\put(35,06){\line(0,1){7}} \put(39,06){\line(0,1){7}}

\put(36,06){\line(0,1){7}} \put(37,06){\line(0,1){7}}
\put(38,06){\line(0,1){7}}

\put(64,06){\line(1,0){04}} \put(64,13){\line(1,0){04}}
\put(64,06){\line(0,1){7}} \put(68,06){\line(0,1){7}}

\put(65,06){\line(0,1){7}} \put(66,06){\line(0,1){7}}
\put(67,06){\line(0,1){7}}


\end{picture}
\end{center}

\section{Combinatorial Engineering Frameworks}

 The suggested combinatorial engineering
 frameworks (as basic ``design frameworks'')
 can be used as support tools
 at various stages of system life cycle (Fig. 6).
 The extended list of the examined
  combinatorial engineering frameworks for  modular systems
 is the following (e.g., \cite{lev06,levprob07}):

 {\bf 1.} Design of a hierarchical system model
 ~(\(T_{1}\)).

 {\bf 2.} Hierarchical modular system design ~(\(T_{2}\)):

  {\it 2.1.} basic hierarchical modular system design
   to obtain a system version
  ~(\(T_{21}\)),

  {\it 2.2.} hierarchical modular system design
 to obtain a family of system versions
  ~(\(T_{22}\)).

  {\bf 3.} Evaluation of system
    (comparison, diagnostics, {\it etc.}) ~(\(T_{3}\)).

 {\bf 4.} Detection of system bottlenecks ~(\(T_{4}\)).
                  \index{detection of system bottlenecks}

 {\bf 5.} Redesign (improvement, upgrade, adaptation)
 ~(\(T_{5}\)):

 {\it 5.1.} basic system improvement
  (``1-1'') ~(\(T_{51}\)),

 {\it 5.2.} system improvement to obtain a family of system versions
 (``1-m'') ~(\(T_{52}\)),

 {\it 5.3.} basic aggregation of system versions into a resultant
 (aggregated) system (``n-1'') ~(\(T_{53}\)),

 {\it 5.4.} aggregation of system versions into a resultant
 (aggregated) system (``n-m'') ~(\(T_{54}\)).

 {\bf 6.} Multistage design
 (i.e., design of a system trajectory)
 ~(\(T_{6}\)).

 {\bf 7.} Modeling of system development/evolution process
 (flow of system \index{system evolution}
  generations) and forecasting ~(\(T_{7}\)).

~~

~~

~~


\begin{center}
\begin{picture}(108,169)

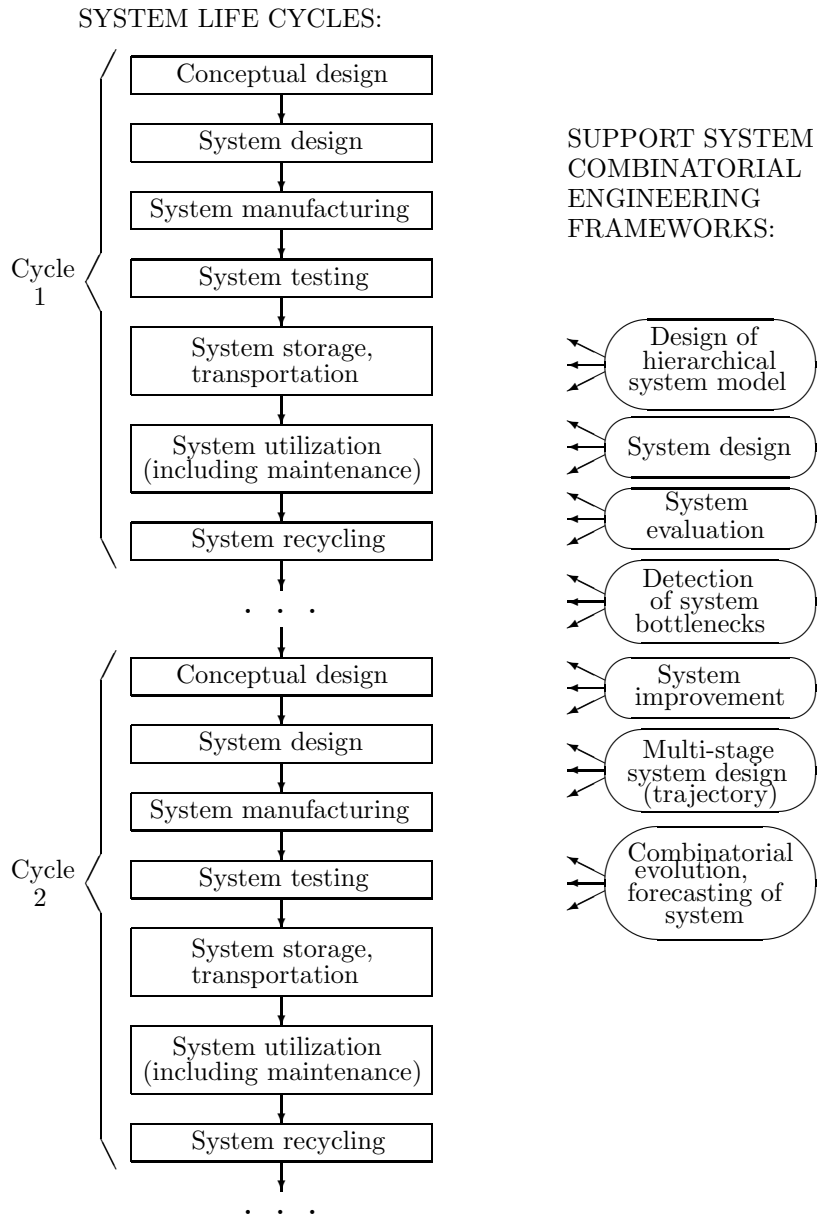
\put(05,00){\makebox(0,0)[bl]{Fig. 6. Life cycles and
combinatorial engineering  frameworks}}


\put(74,148){\makebox(0,0)[bl]{SUPPORT SYSTEM}}
\put(74,144){\makebox(0,0)[bl]{COMBINATORIAL}}
\put(74,140){\makebox(0,0)[bl]{ENGINEERING}}
\put(74,136){\makebox(0,0)[bl]{FRAMEWORKS:}}


\put(93,119){\oval(28,12)}

\put(85,121){\makebox(0,0)[bl]{Design of}}
\put(84,118.4){\makebox(0,0)[bl]{hierarchical}}
\put(82.2,115){\makebox(0,0)[bl]{system model}}

\put(79,119){\vector(-1,0){5}}

\put(79,120){\vector(-2,1){5}} \put(79,118){\vector(-2,-1){5}}


\put(93,108){\oval(28,08)}

\put(82,106.3){\makebox(0,0)[bl]{System design}}

\put(79,108){\vector(-1,0){5}}

\put(79,109){\vector(-2,1){5}} \put(79,107){\vector(-2,-1){5}}


\put(93,98.5){\oval(28,08)}

\put(87,99){\makebox(0,0)[bl]{System}}
\put(84.5,96){\makebox(0,0)[bl]{evaluation}}

\put(79,98.5){\vector(-1,0){5}}

\put(79,99.5){\vector(-2,1){5}} \put(79,97.5){\vector(-2,-1){5}}


\put(93,87.5){\oval(28,11)}

\put(84,89.5){\makebox(0,0)[bl]{Detection}}
\put(85,85.8){\makebox(0,0)[bl]{of system}}
\put(83,83.5){\makebox(0,0)[bl]{bottlenecks}}

\put(79,87.5){\vector(-1,0){5}}

\put(79,88.5){\vector(-2,1){5}} \put(79,86.5){\vector(-2,-1){5}}


\put(93,76){\oval(28,08)}

\put(86,76){\makebox(0,0)[bl]{System}}
\put(83,73){\makebox(0,0)[bl]{improvement}}

\put(79,76){\vector(-1,0){5}}

\put(79,77){\vector(-2,1){5}} \put(79,75){\vector(-2,-1){5}}


\put(93,65){\oval(28,11)}

\put(84,66){\makebox(0,0)[bl]{Multi-stage}}
\put(82,63){\makebox(0,0)[bl]{system design}}
\put(84,60){\makebox(0,0)[bl]{(trajectory)}}

\put(79,65){\vector(-1,0){5}}

\put(79,66){\vector(-2,1){5}} \put(79,64){\vector(-2,-1){5}}


\put(93,50){\oval(28,15)}

\put(82,53){\makebox(0,0)[bl]{Combinatorial}}
\put(83,50){\makebox(0,0)[bl]{evolution,}}
\put(82,47){\makebox(0,0)[bl]{forecasting of}}
\put(87,44){\makebox(0,0)[bl]{system }}

\put(79,50){\vector(-1,0){5}}

\put(79,51){\vector(-2,1){5}} \put(79,49){\vector(-2,-1){5}}


\put(09,164){\makebox(0,0)[bl]{SYSTEM LIFE CYCLES: }}


\put(16,155){\line(1,0){40}}\put(16,160){\line(1,0){40}}
\put(16,155){\line(0,1){05}} \put(56,155){\line(0,1){05}}

\put(22,156){\makebox(0,0)[bl]{Conceptual design}}

\put(36,155){\vector(0,-1){4}}


\put(16,146){\line(1,0){40}}\put(16,151){\line(1,0){40}}
\put(16,146){\line(0,1){05}} \put(56,146){\line(0,1){05}}

\put(25,147){\makebox(0,0)[bl]{System design}}

\put(36,146){\vector(0,-1){4}}


\put(16,137){\line(1,0){40}}\put(16,142){\line(1,0){40}}
\put(16,137){\line(0,1){05}} \put(56,137){\line(0,1){05}}

\put(18.5,138){\makebox(0,0)[bl]{System manufacturing}}

\put(36,137){\vector(0,-1){4}}


\put(16,128){\line(1,0){40}}\put(16,133){\line(1,0){40}}
\put(16,128){\line(0,1){05}} \put(56,128){\line(0,1){05}}

\put(25,129){\makebox(0,0)[bl]{System testing}}

\put(36,128){\vector(0,-1){4}}


\put(00,130){\makebox(0,0)[bl]{Cycle}}
\put(03,127){\makebox(0,0)[bl]{\(1\)}}

\put(10,130){\line(1,2){2}} \put(10,130){\line(1,-2){2}}

\put(12,134){\line(0,1){23}} \put(12,157){\line(1,2){2}}

\put(12,126){\line(0,-1){30}} \put(12,96){\line(1,-2){2}}



\put(16,115){\line(1,0){40}}\put(16,124){\line(1,0){40}}
\put(16,115){\line(0,1){09}} \put(56,115){\line(0,1){09}}

\put(24,119.5){\makebox(0,0)[bl]{System storage,}}
\put(24,116){\makebox(0,0)[bl]{transportation}}

\put(36,115){\vector(0,-1){4}}


\put(16,102){\line(1,0){40}}\put(16,111){\line(1,0){40}}
\put(16,102){\line(0,1){09}} \put(56,102){\line(0,1){09}}

\put(21.5,106.5){\makebox(0,0)[bl]{System utilization}}
\put(17.5,103){\makebox(0,0)[bl]{(including maintenance)}}

\put(36,102){\vector(0,-1){4}}


\put(16,93){\line(1,0){40}}\put(16,98){\line(1,0){40}}
\put(16,93){\line(0,1){05}} \put(56,93){\line(0,1){05}}

\put(24,94){\makebox(0,0)[bl]{System recycling}}

\put(36,93){\vector(0,-1){4}}


\put(31,86){\makebox(0,0)[bl]{{{\bf .~ .~ .}}}}

\put(36,84){\vector(0,-1){4}}


\put(16,75){\line(1,0){40}}\put(16,80){\line(1,0){40}}
\put(16,75){\line(0,1){05}} \put(56,75){\line(0,1){05}}

\put(22,76){\makebox(0,0)[bl]{Conceptual design}}

\put(36,75){\vector(0,-1){4}}


\put(16,66){\line(1,0){40}}\put(16,71){\line(1,0){40}}
\put(16,66){\line(0,1){05}} \put(56,66){\line(0,1){05}}

\put(25,67){\makebox(0,0)[bl]{System design}}

\put(36,66){\vector(0,-1){4}}


\put(16,57){\line(1,0){40}}\put(16,62){\line(1,0){40}}
\put(16,57){\line(0,1){05}} \put(56,57){\line(0,1){05}}

\put(18.5,58){\makebox(0,0)[bl]{System manufacturing}}

\put(36,57){\vector(0,-1){4}}


\put(16,48){\line(1,0){40}}\put(16,53){\line(1,0){40}}
\put(16,48){\line(0,1){05}} \put(56,48){\line(0,1){05}}

\put(25,49){\makebox(0,0)[bl]{System testing}}

\put(36,48){\vector(0,-1){4}}


\put(0,50){\makebox(0,0)[bl]{Cycle}}
\put(03,47){\makebox(0,0)[bl]{\(2\)}}

\put(10,50){\line(1,2){2}} \put(10,50){\line(1,-2){2}}

\put(12,54){\line(0,1){23}} \put(12,77){\line(1,2){2}}

\put(12,46){\line(0,-1){30}} \put(12,16){\line(1,-2){2}}



\put(16,35){\line(1,0){40}}\put(16,44){\line(1,0){40}}
\put(16,35){\line(0,1){09}} \put(56,35){\line(0,1){09}}

\put(24,39.5){\makebox(0,0)[bl]{System storage,}}
\put(24,36){\makebox(0,0)[bl]{transportation}}

\put(36,35){\vector(0,-1){4}}


\put(16,22){\line(1,0){40}}\put(16,31){\line(1,0){40}}
\put(16,22){\line(0,1){09}} \put(56,22){\line(0,1){09}}

\put(21.5,26.5){\makebox(0,0)[bl]{System utilization}}
\put(17.5,23){\makebox(0,0)[bl]{(including maintenance)}}

\put(36,22){\vector(0,-1){4}}


\put(16,13){\line(1,0){40}}\put(16,18){\line(1,0){40}}
\put(16,13){\line(0,1){05}} \put(56,13){\line(0,1){05}}

\put(24,14){\makebox(0,0)[bl]{System recycling}}

\put(36,13){\vector(0,-1){4}}


\put(31,06){\makebox(0,0)[bl]{{{\bf .~ .~ .}}}}

\end{picture}
\end{center}

  The frameworks above can be applied to systems, systems
 requirements, standards, plans, etc.
 (e.g., \cite{lev98,lev06,lev11agg}).
 A generalized scheme of our research domain is presented in Fig. 7.

  Mainly, several combinatorial engineering
 frameworks are often used together in applications,
 for example:

 (i) design of system hierarchical model, system design,
 detection of system bottlenecks,
 system improvement;

 (ii) design of system hierarchical model,
 detection of system bottlenecks,
  combinatorial evolution of the system,
 design of system forecasts, aggregation of the forecasts.

\begin{center}
\begin{picture}(117,116)

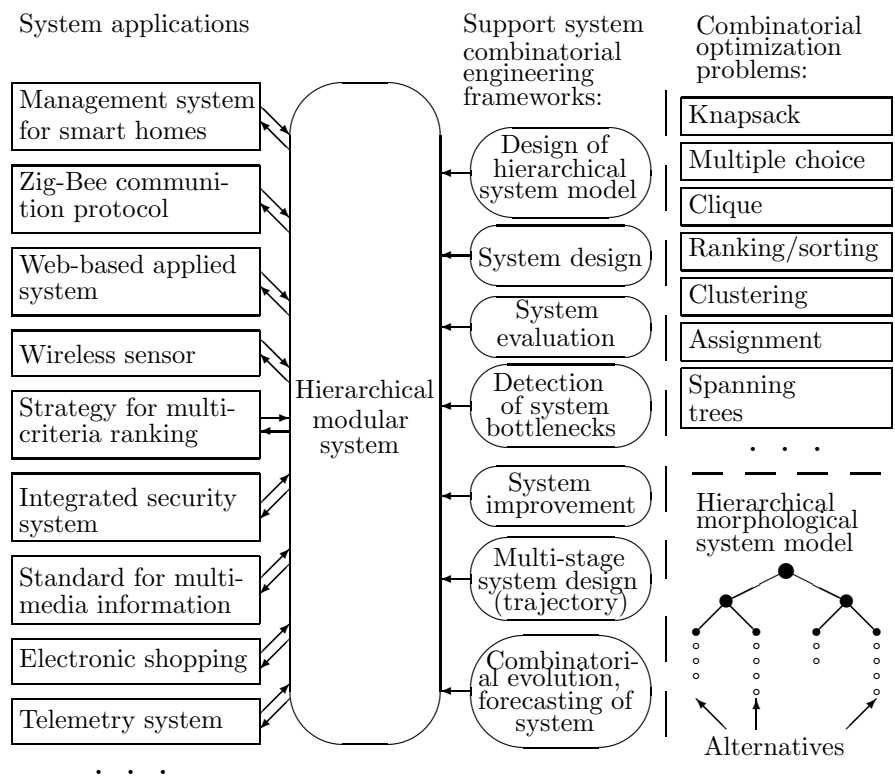
\put(17,00){\makebox(0,0)[bl]{Fig. 7. Generalized scheme of
 examined domain}}


\put(91,107.5){\makebox(0,0)[bl]{Combinatorial}}
\put(91,104){\makebox(0,0)[bl]{optimization}}
\put(91,101){\makebox(0,0)[bl]{problems:}}


\put(89,94){\line(1,0){28}} \put(89,99){\line(1,0){28}}
\put(89,94){\line(0,1){05}} \put(117,94){\line(0,1){05}}

\put(90,95){\makebox(0,0)[bl]{Knapsack}}


\put(89,88){\line(1,0){28}} \put(89,93){\line(1,0){28}}
\put(89,88){\line(0,1){5}} \put(117,88){\line(0,1){05}}

\put(90,89){\makebox(0,0)[bl]{Multiple choice}}


\put(89,82){\line(1,0){28}} \put(89,87){\line(1,0){28}}
\put(89,82){\line(0,1){5}} \put(117,82){\line(0,1){05}}

\put(90,83){\makebox(0,0)[bl]{Clique}}


\put(89,76){\line(1,0){28}} \put(89,81){\line(1,0){28}}
\put(89,76){\line(0,1){5}} \put(117,76){\line(0,1){05}}

\put(90,77){\makebox(0,0)[bl]{Ranking/sorting}}


\put(89,70){\line(1,0){28}} \put(89,75){\line(1,0){28}}
\put(89,70){\line(0,1){5}} \put(117,70){\line(0,1){05}}

\put(90,71){\makebox(0,0)[bl]{Clustering}}


\put(89,64){\line(1,0){28}} \put(89,69){\line(1,0){28}}
\put(89,64){\line(0,1){5}} \put(117,64){\line(0,1){05}}

\put(90,65){\makebox(0,0)[bl]{Assignment}}


\put(89,55){\line(1,0){28}} \put(89,63){\line(1,0){28}}
\put(89,55){\line(0,1){08}} \put(117,55){\line(0,1){08}}

\put(90,59){\makebox(0,0)[bl]{Spanning }}
\put(90,56){\makebox(0,0)[bl]{trees }}

\put(98,52){\makebox(0,0)[bl]{{\bf .~ .~ .}}}


\put(90.5,49){\line(1,0){4}} \put(97.5,49){\line(1,0){4}}
\put(104.5,49){\line(1,0){4}} \put(111.5,49){\line(1,0){4}}

\put(91,44){\makebox(0,0)[bl]{Hierarchical}}
\put(91,41){\makebox(0,0)[bl]{morphological}}
\put(91,38){\makebox(0,0)[bl]{system model}}

\put(103,36){\circle*{2.0}}

\put(103,36){\line(-2,-1){8}} \put(103,36){\line(2,-1){8}}

\put(95,32){\circle*{1.6}} \put(111,32){\circle*{1.6}}

\put(95,32){\line(-1,-1){4}} \put(95,32){\line(1,-1){4}}

\put(91,28){\circle*{1.0}} \put(99,28){\circle*{1.0}}

\put(111,32){\line(-1,-1){4}} \put(111,32){\line(1,-1){4}}

\put(107,28){\circle*{1.0}} \put(115,28){\circle*{1.0}}


\put(91,26){\circle{0.8}} \put(91,24){\circle{0.8}}
\put(91,22){\circle{0.8}}


\put(99,26){\circle{0.8}} \put(99,24){\circle{0.8}}
\put(99,22){\circle{0.8}} \put(99,20){\circle{0.8}}

\put(107,26){\circle{0.8}} \put(107,24){\circle{0.8}}

\put(115,26){\circle{0.8}} \put(115,24){\circle{0.8}}
\put(115,22){\circle{0.8}} \put(115,20){\circle{0.8}}

\put(92,11.6){\makebox(0,0)[bl]{Alternatives}}

\put(95,15){\vector(-1,1){4}} \put(99,15){\vector(0,1){4}}

\put(111,15){\vector(1,1){4}}


\put(87,100){\line(0,-1){6}} \put(87,90){\line(0,-1){6}}
\put(87,80){\line(0,-1){6}} \put(87,70){\line(0,-1){6}}
\put(87,60){\line(0,-1){6}} \put(87,50){\line(0,-1){6}}
\put(87,40){\line(0,-1){5}} \put(87,30){\line(0,-1){6}}
\put(87,20){\line(0,-1){6}}


\put(60,107){\makebox(0,0)[bl]{Support system}}
\put(60,104){\makebox(0,0)[bl]{combinatorial}}
\put(60,100.5){\makebox(0,0)[bl]{engineering}}
\put(60,98){\makebox(0,0)[bl]{frameworks:}}


\put(73,89){\oval(24,12)}

\put(65,91){\makebox(0,0)[bl]{Design of}}
\put(64,88.4){\makebox(0,0)[bl]{hierarchical}}
\put(62.2,85){\makebox(0,0)[bl]{system model}}

\put(61,89){\vector(-1,0){4}}


\put(73,78){\oval(24,08)}

\put(62,76){\makebox(0,0)[bl]{System design}}

\put(61,78){\vector(-1,0){4}}


\put(73,68.5){\oval(24,08)}

\put(67,69){\makebox(0,0)[bl]{System}}
\put(64.5,66){\makebox(0,0)[bl]{evaluation}}

\put(61,68.5){\vector(-1,0){4}}


\put(73,58){\oval(24,11)}

\put(64,60){\makebox(0,0)[bl]{Detection}}
\put(65,56.3){\makebox(0,0)[bl]{of system}}
\put(63,54){\makebox(0,0)[bl]{bottlenecks}}

\put(61,58){\vector(-1,0){4}}


\put(73,46){\oval(24,08)}

\put(66,46){\makebox(0,0)[bl]{System}}
\put(63,43){\makebox(0,0)[bl]{improvement}}

\put(61,46){\vector(-1,0){4}}


\put(73,35){\oval(24,11)}

\put(64,36){\makebox(0,0)[bl]{Multi-stage}}
\put(62,33){\makebox(0,0)[bl]{system design}}
\put(64,30){\makebox(0,0)[bl]{(trajectory)}}

\put(61,35){\vector(-1,0){4}}


\put(73,20){\oval(24,15)}

\put(63,23){\makebox(0,0)[bl]{Combinatori-}}
\put(62,20){\makebox(0,0)[bl]{al evolution,}}
\put(62,17){\makebox(0,0)[bl]{forecasting of}}
\put(67,14){\makebox(0,0)[bl]{system }}

\put(61,20){\vector(-1,0){4}}


\put(47,57){\oval(20,88)}

\put(37.7,59){\makebox(0,0)[bl]{Hierarchical}}
\put(40,55){\makebox(0,0)[bl]{modular}}
\put(41,51){\makebox(0,0)[bl]{system}}


\put(01,107){\makebox(0,0)[bl]{System applications}}


\put(00,92){\line(1,0){33}}\put(00,101){\line(1,0){33}}
\put(00,92){\line(0,1){09}} \put(33,92){\line(0,1){09}}

\put(01,97){\makebox(0,0)[bl]{Management system}}
\put(01,93.5){\makebox(0,0)[bl]{for smart homes}}

\put(33,98){\vector(1,-1){4}} \put(37,92){\vector(-1,1){4}}


\put(00,81){\line(1,0){33}}\put(00,90){\line(1,0){33}}
\put(00,81){\line(0,1){09}} \put(33,81){\line(0,1){09}}

\put(01,86){\makebox(0,0)[bl]{Zig-Bee communi-}}
\put(01,82.5){\makebox(0,0)[bl]{tion protocol}}

\put(33,87){\vector(1,-1){4}} \put(37,81){\vector(-1,1){4}}


\put(00,70){\line(1,0){33}} \put(00,79){\line(1,0){33}}
\put(00,70){\line(0,1){09}} \put(33,70){\line(0,1){09}}

\put(01,75){\makebox(0,0)[bl]{Web-based applied}}
\put(01,71.5){\makebox(0,0)[bl]{system}}

\put(33,76){\vector(1,-1){4}} \put(37,70){\vector(-1,1){4}}


\put(00,62){\line(1,0){33}}\put(00,68){\line(1,0){33}}
\put(00,62){\line(0,1){06}} \put(33,62){\line(0,1){06}}

\put(01,63.7){\makebox(0,0)[bl]{Wireless sensor}}

\put(33,67){\vector(1,-1){4}} \put(37,61){\vector(-1,1){4}}


\put(00,51){\line(1,0){33}}\put(00,60){\line(1,0){33}}
\put(00,51){\line(0,1){09}} \put(33,51){\line(0,1){09}}

\put(01,56){\makebox(0,0)[bl]{Strategy for multi-}}
\put(01,52.5){\makebox(0,0)[bl]{criteria ranking}}

\put(33,56.4){\vector(1,0){4}} \put(37,54.6){\vector(-1,0){4}}


\put(00,40){\line(1,0){33}}\put(00,49){\line(1,0){33}}
\put(00,40){\line(0,1){09}} \put(33,40){\line(0,1){09}}

\put(01,44){\makebox(0,0)[bl]{Integrated security}}
\put(01,40.5){\makebox(0,0)[bl]{system}}

\put(33,45){\vector(1,1){4}} \put(37,47){\vector(-1,-1){4}}


\put(00,29){\line(1,0){33}}\put(00,38){\line(1,0){33}}
\put(00,29){\line(0,1){09}} \put(33,29){\line(0,1){09}}

\put(01,34){\makebox(0,0)[bl]{Standard for multi-}}
\put(01,30.5){\makebox(0,0)[bl]{media information}}

\put(33,35){\vector(1,1){4}} \put(37,37){\vector(-1,-1){4}}


\put(00,21){\line(1,0){33}}\put(00,27){\line(1,0){33}}
\put(00,21){\line(0,1){06}} \put(33,21){\line(0,1){06}}

\put(01,22.5){\makebox(0,0)[bl]{Electronic shopping}}

\put(33,25){\vector(1,1){4}} \put(37,27){\vector(-1,-1){4}}


\put(00,13){\line(1,0){33}}\put(00,19){\line(1,0){33}}
\put(00,13){\line(0,1){06}} \put(33,13){\line(0,1){06}}

\put(01,14.5){\makebox(0,0)[bl]{Telemetry system}}

\put(33,17){\vector(1,1){4}} \put(37,19){\vector(-1,-1){4}}


\put(11,09){\makebox(0,0)[bl]{{{\bf .~ .~ .}}}}

\end{picture}
\end{center}

\section{Conclusion}

 This paper describes a methodological viewpoint
 to the set of basic system combinatorial engineering frameworks
 for maintenance (modeling, design, improvement, etc.)
 of modular systems with a hierarchical morphological structure.
 Table 2 contains
 the list of combinatorial engineering
 frameworks and corresponding underlying combinatorial optimization problems
 and combinatorial engineering frameworks.
  It is necessary to note
  the  system family design/improvement frameworks
 (\(T_{22}\), \(T_{52}\), \(T_{54}\))
  require special additional research efforts.
%
%
 Some author's system design applications
 based on the combinatorial engineering frameworks are pointed out
 in Table 3.
 Evidently, the considered combinatorial engineering
 frameworks can be successfully applied in
 education
 (engineering, management, computer science, applied mathematics)
 including student-project based courses in system design
 (e.g., \cite{lev06a,lev09edu,lev11ed,lev11csedu}).

%
 In the future, it may be prospective
  to consider the following research directions:

 {\bf I.} Methodology:

 {\it 1.1.} examination of various network-like models
 (e.g., Petri nets)
  for modular systems;

 {\it 1.2.} further investigation of system evolution/development
 processes and system forecasting;

  {\it 1.3.} suggestion and investigation of
 combinatorial engineering frameworks
 (as special macro-heuristics)
 for well-known complex combinatorial optimization
 problems (e.g., timetabling, augmentation problem);

 {\it 1.4.} consideration  of uncertainty and  AI techniques;

 {\it 1.5.} special studies of dynamical systems.

 {\bf II.} Applications:

 {\it 2.1.} consideration of various applied systems
 (e.g., power systems, communication systems,
 various applied systems in social engineering);

 {\it 2.2.} special studies have to be targeted to financial
 engineering;

 {\it 2.3.} special tools have to be designed for
 individual modular educational systems for student usage;

 {\it 2.4.} special efforts have to be targeted to
 biomedical applications
 (e.g., diagnosis, medical treatment, etc.).

\begin{center}
\begin{picture}(121,164)

\put(07,160){\makebox(0,0)[bl]{Table 2. Combinatorial
 engineering frameworks and their description}}

\put(00,00){\line(1,0){121}} \put(00,141.5){\line(1,0){121}}
\put(60,152){\line(1,0){61}} \put(00,158){\line(1,0){121}}

\put(00,00){\line(0,1){158}} \put(60,00){\line(0,1){158}}
\put(85,00){\line(0,1){152}} \put(121,00){\line(0,1){158}}

\put(01,154){\makebox(0,0)[bl]{Combinatorial engineering}}
\put(01,150){\makebox(0,0)[bl]{framework}}

\put(80,153.5){\makebox(0,0)[bl]{Description}}

\put(61,147){\makebox(0,0)[bl]{General design}}
\put(61,143.3){\makebox(0,0)[bl]{framework}}


\put(86,147){\makebox(0,0)[bl]{Underlying problems,}}
\put(86,143){\makebox(0,0)[bl]{frameworks}}


\put(01,136){\makebox(0,0)[bl]{1.Design of system hierarchical }}

\put(03,132){\makebox(0,0)[bl]{model
 ~(\(T_{1}\)) }}

\put(61,136){\makebox(0,0)[bl]{\cite{lev12hier}}}

\put(86,136){\makebox(0,0)[bl]{Spanning trees,}}
\put(86,132){\makebox(0,0)[bl]{assignment/allocation,}}
\put(86,128){\makebox(0,0)[bl]{clustering, multilevel}}
\put(86,124){\makebox(0,0)[bl]{network design}}


\put(01,119){\makebox(0,0)[bl]{2.System  design ~(\(T_{2}\)):}}


\put(01,103){\makebox(0,0)[bl]{2.1.basic system design}}

\put(05,99){\makebox(0,0)[bl]{(one resultant version)
  ~(\(T_{21}\))}}

\put(61,119){\makebox(0,0)[bl]{\cite{lev98,lev06,lev09,lev11agg,lev12morph}}}
\put(61,115){\makebox(0,0)[bl]{\cite{lev12a}}}

\put(86,119){\makebox(0,0)[bl]{Morphological clique,}}
\put(86,115){\makebox(0,0)[bl]{multiple choice, }}
\put(86,111){\makebox(0,0)[bl]{assignment/allocation, }}
\put(86,107){\makebox(0,0)[bl]{agreement problems }}


\put(01,95){\makebox(0,0)[bl]{2.2.system family design }}

\put(05,91){\makebox(0,0)[bl]{(several resultant versions)
 ~(\(T_{22}\)) }}


\put(01,86){\makebox(0,0)[bl]{3.System evaluation ~(\(T_{3}\))}}

\put(61,86){\makebox(0,0)[bl]{\cite{lev06}}}

\put(86,86){\makebox(0,0)[bl]{Multicriteria ranking}}


\put(01,81){\makebox(0,0)[bl]{4.Detection of bottlenecks
 ~(\(T_{4}\)) }}

\put(61,81){\makebox(0,0)[bl]{\cite{lev98,lev06}}}

\put(86,81.5){\makebox(0,0)[bl]{Detection of critical}}
\put(86,77){\makebox(0,0)[bl]{components (e.g.,}}
\put(86,73){\makebox(0,0)[bl]{multicriteria ranking,}}
\put(86,69){\makebox(0,0)[bl]{dominating set)}}


\put(01,64){\makebox(0,0)[bl]{5.System improvement
 ~(\(T_{5}\)):}}

\put(61,64){\makebox(0,0)[bl]{\cite{lev98,lev06,lev09,levsib10,lev11agg}}}
\put(61,60){\makebox(0,0)[bl]{\cite{lev12morph,lev12a}}}

\put(86,64){\makebox(0,0)[bl]{\(T_{1},T_{2},T_{3},T_{4}\)}}


\put(01,56){\makebox(0,0)[bl]{5.1.basic system improvement,}}

\put(05,52){\makebox(0,0)[bl]{result:  one version (``1-1'')
  ~(\(T_{51}\)) }}


\put(01,48){\makebox(0,0)[bl]{5.2.system improvement, result: }}

\put(05,44){\makebox(0,0)[bl]{several system versions }}

\put(05,40){\makebox(0,0)[bl]{(``1-m'') ~(\(T_{52}\)) }}


\put(01,36){\makebox(0,0)[bl]{5.3.aggregation of system
versions:}}

\put(05,32){\makebox(0,0)[bl]{one resultant (aggregated)}}

\put(05,28){\makebox(0,0)[bl]{system (``n-1'')
  ~(\(T_{53}\)) }}


\put(01,24){\makebox(0,0)[bl]{5.4.aggregation of system
versions:}}

\put(05,20){\makebox(0,0)[bl]{several resultant (aggregated)}}

\put(05,16){\makebox(0,0)[bl]{systems (``n-m'')
 ~(\(T_{54}\)) }}


\put(01,11){\makebox(0,0)[bl]{6.Multistage system design }}

\put(05,07){\makebox(0,0)[bl]{(design of system trajectory)
 ~(\(T_{6}\))}}

\put(61,11){\makebox(0,0)[bl]{\cite{lev98,lev06,lev09,lev11agg,lev12a}}}

\put(86,11){\makebox(0,0)[bl]{\(T_{1},T_{2}\)}}


\put(01,02){\makebox(0,0)[bl]{7.System evolution, forecasting
 ~(\(T_{7}\) )}}

\put(61,02){\makebox(0,0)[bl]{\cite{lev06,lev13}}}

\put(86,02){\makebox(0,0)[bl]{\(T_{1},T_{2},T_{4},T_{5}\)}}


\end{picture}
\end{center}

\begin{center}
\begin{picture}(115,125)
\put(02,121){\makebox(0,0)[bl]{Table 3. System applications and
combinatorial engineering frameworks}}

\put(1,114.5){\makebox(0,0)[bl]{System application}}

\put(52.5,115){\makebox(0,0)[bl]{Used frameworks}}

\put(00,00){\line(1,0){115}} \put(00,107){\line(1,0){115}}
\put(45,113){\line(1,0){42}} \put(00,119){\line(1,0){115}}

\put(00,00){\line(0,1){119}} \put(45,00){\line(0,1){119}}
\put(87,00){\line(0,1){119}} \put(115,00){\line(0,1){119}}

\put(51,107){\line(0,1){6}} \put(57,107){\line(0,1){6}}
\put(63,107){\line(0,1){6}} \put(69,107){\line(0,1){6}}
\put(75,107){\line(0,1){6}} \put(81,107){\line(0,1){6}}

\put(46,108.5){\makebox(0,0)[bl]{\(T_{1}\)}}
\put(52,108.5){\makebox(0,0)[bl]{\(T_{2}\)}}
\put(58,108.5){\makebox(0,0)[bl]{\(T_{3}\)}}
\put(64,108.5){\makebox(0,0)[bl]{\(T_{4}\)}}
\put(70,108.5){\makebox(0,0)[bl]{\(T_{5}\)}}
\put(76,108.5){\makebox(0,0)[bl]{\(T_{6}\)}}
\put(82,108.5){\makebox(0,0)[bl]{\(T_{7}\)}}

\put(88,115){\makebox(0,0)[bl]{Source}}

\put(01,102){\makebox(0,0)[bl]{1.Electronic shopping}}

\put(47,102){\makebox(0,0)[bl]{\(\star \)}}
\put(53,102){\makebox(0,0)[bl]{\(\star \)}}
\put(59,102){\makebox(0,0)[bl]{\(\star \)}}
\put(65,102){\makebox(0,0)[bl]{\(\star \)}}
\put(71,102){\makebox(0,0)[bl]{\(\star \)}}
\put(77,102){\makebox(0,0)[bl]{\(\star \)}}
\put(83,102){\makebox(0,0)[bl]{}}

\put(88,101.5){\makebox(0,0)[bl]{\cite{lev12shop}}}

\put(01,98){\makebox(0,0)[bl]{2.Web-based
 system}}

\put(47,98){\makebox(0,0)[bl]{\(\star \)}}
\put(53,98){\makebox(0,0)[bl]{\(\star \)}}
\put(59,98){\makebox(0,0)[bl]{\(\star \)}}
\put(65,98){\makebox(0,0)[bl]{\(\star \)}}
\put(71,98){\makebox(0,0)[bl]{\(\star \)}}
\put(77,98){\makebox(0,0)[bl]{\(\star \)}}
\put(83,98){\makebox(0,0)[bl]{}}

\put(88,97.5){\makebox(0,0)[bl]{\cite{lev11inf,lev11agg}}}

\put(01,94){\makebox(0,0)[bl]{3.Strategy for sorting}}
\put(04,90){\makebox(0,0)[bl]{(multicriteria ranking)}}

\put(47,94){\makebox(0,0)[bl]{\(\star \)}}
\put(53,94){\makebox(0,0)[bl]{\(\star \)}}
\put(59,94){\makebox(0,0)[bl]{}} \put(65,90){\makebox(0,0)[bl]{}}
\put(71,94){\makebox(0,0)[bl]{}} \put(77,90){\makebox(0,0)[bl]{}}
\put(83,94){\makebox(0,0)[bl]{}}

\put(88,93.5){\makebox(0,0)[bl]{\cite{lev98,lev12b}}}

\put(01,86.5){\makebox(0,0)[bl]{4.DSS COMBI}}

\put(47,86){\makebox(0,0)[bl]{\(\star \)}}
\put(53,86){\makebox(0,0)[bl]{}} \put(59,78){\makebox(0,0)[bl]{}}
\put(65,86){\makebox(0,0)[bl]{}} \put(71,78){\makebox(0,0)[bl]{}}
\put(77,86){\makebox(0,0)[bl]{}}
\put(83,86){\makebox(0,0)[bl]{\(\star \)}}

\put(88,85.5){\makebox(0,0)[bl]{\cite{lev93,lev98}}}

\put(01,82.5){\makebox(0,0)[bl]{5.Modular software}}

\put(47,82){\makebox(0,0)[bl]{\(\star \)}}
\put(53,82){\makebox(0,0)[bl]{\(\star \)}}
\put(59,82){\makebox(0,0)[bl]{\(\star \)}}
\put(65,82){\makebox(0,0)[bl]{\(\star \)}}
\put(71,82){\makebox(0,0)[bl]{\(\star \)}}
\put(77,82){\makebox(0,0)[bl]{}} \put(83,78){\makebox(0,0)[bl]{}}

\put(88,81.5){\makebox(0,0)[bl]{\cite{lev05,lev06}}}

\put(01,78.5){\makebox(0,0)[bl]{6.Notebook }}

\put(47,78){\makebox(0,0)[bl]{\(\star \)}}
\put(53,78){\makebox(0,0)[bl]{\(\star \)}}
\put(59,78){\makebox(0,0)[bl]{}} \put(65,78){\makebox(0,0)[bl]{}}
\put(71,78){\makebox(0,0)[bl]{\(\star \)}}
\put(77,78){\makebox(0,0)[bl]{}} \put(83,78){\makebox(0,0)[bl]{}}

\put(88,77.5){\makebox(0,0)[bl]{\cite{lev11agg}}}

\put(01,74){\makebox(0,0)[bl]{7.Regional network}}

\put(47,74){\makebox(0,0)[bl]{\(\star \)}}
\put(53,74){\makebox(0,0)[bl]{}}
\put(59,74){\makebox(0,0)[bl]{\(\star \)}}
\put(65,74){\makebox(0,0)[bl]{}} \put(71,74){\makebox(0,0)[bl]{}}
\put(77,74){\makebox(0,0)[bl]{}} \put(83,74){\makebox(0,0)[bl]{}}

\put(88,73.5){\makebox(0,0)[bl]{\cite{levsaf11}}}

\put(01,70.5){\makebox(0,0)[bl]{8.GSM network }}

\put(47,70){\makebox(0,0)[bl]{\(\star \)}}
\put(53,70){\makebox(0,0)[bl]{\(\star \)}}
\put(59,70){\makebox(0,0)[bl]{}} \put(65,70){\makebox(0,0)[bl]{}}
\put(71,70){\makebox(0,0)[bl]{}} \put(77,70){\makebox(0,0)[bl]{}}
\put(83,70){\makebox(0,0)[bl]{}}

\put(88,69.5){\makebox(0,0)[bl]{\cite{lev12morph,levvis07}}}

\put(01,66){\makebox(0,0)[bl]{9.Telemetry system }}

\put(47,66){\makebox(0,0)[bl]{\(\star \)}}
\put(53,66){\makebox(0,0)[bl]{\(\star \)}}
\put(59,66){\makebox(0,0)[bl]{\(\star \)}}
\put(65,66){\makebox(0,0)[bl]{\(\star \)}}
\put(71,66){\makebox(0,0)[bl]{\(\star \)}}
\put(77,66){\makebox(0,0)[bl]{}} \put(83,66){\makebox(0,0)[bl]{}}

\put(88,65.5){\makebox(0,0)[bl]{\cite{lev12tele,levkhod07}}}

\put(0.5,62){\makebox(0,0)[bl]{10.Security system}}

\put(47,62){\makebox(0,0)[bl]{\(\star \)}}
\put(53,62){\makebox(0,0)[bl]{\(\star \)}}
\put(59,62){\makebox(0,0)[bl]{\(\star \)}}
\put(65,62){\makebox(0,0)[bl]{\(\star \)}}
\put(71,62){\makebox(0,0)[bl]{\(\star \)}}
\put(77,62){\makebox(0,0)[bl]{}} \put(83,62){\makebox(0,0)[bl]{}}

\put(88,61.5){\makebox(0,0)[bl]{\cite{lev11agg,levleus09}}}

\put(0.5,58.5){\makebox(0,0)[bl]{11.Wireless sensor}}

\put(47,58){\makebox(0,0)[bl]{\(\star \)}}
\put(53,58){\makebox(0,0)[bl]{\(\star \)}}
\put(59,58){\makebox(0,0)[bl]{\(\star \)}}
\put(65,58){\makebox(0,0)[bl]{\(\star \)}}
\put(71,58){\makebox(0,0)[bl]{\(\star \)}}
\put(77,58){\makebox(0,0)[bl]{}} \put(83,58){\makebox(0,0)[bl]{}}

\put(88,57.5){\makebox(0,0)[bl]{\cite{lev11agg,levfim10,levfim12}}}

\put(0.5,54){\makebox(0,0)[bl]{12.System tests}}

\put(47,54){\makebox(0,0)[bl]{\(\star \)}}
\put(53,54){\makebox(0,0)[bl]{\(\star \)}}
\put(59,54){\makebox(0,0)[bl]{}} \put(65,54){\makebox(0,0)[bl]{}}
\put(71,54){\makebox(0,0)[bl]{}} \put(77,54){\makebox(0,0)[bl]{}}
\put(83,54){\makebox(0,0)[bl]{}}

\put(88,53.5){\makebox(0,0)[bl]{\cite{levlast06}}}

\put(0.5,49.6){\makebox(0,0)[bl]{13.Communication
  protocol}}

\put(47,50){\makebox(0,0)[bl]{\(\star \)}}
\put(53,50){\makebox(0,0)[bl]{\(\star \)}}
\put(59,50){\makebox(0,0)[bl]{\(\star \)}}
\put(65,50){\makebox(0,0)[bl]{\(\star \)}}
\put(71,50){\makebox(0,0)[bl]{\(\star \)}}
\put(77,50){\makebox(0,0)[bl]{\(\star \)}}
\put(83,50){\makebox(0,0)[bl]{\(\star \)}}

\put(88,49.5){\makebox(0,0)[bl]{\cite{lev13mpeg,levand10,levand12}}}

\put(0.5,46){\makebox(0,0)[bl]{14.User interface}}

\put(47,46){\makebox(0,0)[bl]{\(\star \)}}
\put(53,46){\makebox(0,0)[bl]{\(\star \)}}
\put(59,46){\makebox(0,0)[bl]{}} \put(65,46){\makebox(0,0)[bl]{}}
\put(71,46){\makebox(0,0)[bl]{}} \put(77,46){\makebox(0,0)[bl]{}}
\put(83,46){\makebox(0,0)[bl]{}}

\put(88,45.5){\makebox(0,0)[bl]{\cite{lev94,lev98,lev02}}}

\put(0.5,41.5){\makebox(0,0)[bl]{15.Two-floor building}}

\put(47,42){\makebox(0,0)[bl]{\(\star \)}}
\put(53,42){\makebox(0,0)[bl]{\(\star \)}}
\put(59,42){\makebox(0,0)[bl]{\(\star \)}}
\put(65,42){\makebox(0,0)[bl]{\(\star \)}}
\put(71,42){\makebox(0,0)[bl]{}} \put(77,42){\makebox(0,0)[bl]{}}
\put(83,42){\makebox(0,0)[bl]{}}

\put(88,41.5){\makebox(0,0)[bl]{\cite{lev06,levdan05}}}

\put(0.5,38){\makebox(0,0)[bl]{16.Control
 in smart home }}

\put(47,38){\makebox(0,0)[bl]{\(\star \)}}
\put(53,38){\makebox(0,0)[bl]{\(\star \)}}
\put(59,38){\makebox(0,0)[bl]{\(\star \)}}
\put(65,38){\makebox(0,0)[bl]{\(\star \)}}
\put(71,38){\makebox(0,0)[bl]{\(\star \)}}
\put(77,38){\makebox(0,0)[bl]{}} \put(83,38){\makebox(0,0)[bl]{}}

\put(88,37.5){\makebox(0,0)[bl]{\cite{lev13home,levand11a,levand11}}}

\put(0.5,34){\makebox(0,0)[bl]{17.Combinatorial investment}}

\put(47,34){\makebox(0,0)[bl]{\(\star \)}}
\put(53,34){\makebox(0,0)[bl]{\(\star \)}}
\put(59,34){\makebox(0,0)[bl]{}} \put(65,34){\makebox(0,0)[bl]{}}
\put(71,34){\makebox(0,0)[bl]{}} \put(77,34){\makebox(0,0)[bl]{}}
\put(83,34){\makebox(0,0)[bl]{}}

\put(88,33.5){\makebox(0,0)[bl]{\cite{lev98,lev11agg}}}

\put(0.5,29.5){\makebox(0,0)[bl]{18.Political management}}

\put(47,30){\makebox(0,0)[bl]{\(\star \)}}
\put(53,30){\makebox(0,0)[bl]{}} \put(59,30){\makebox(0,0)[bl]{}}
\put(65,30){\makebox(0,0)[bl]{}} \put(71,30){\makebox(0,0)[bl]{}}
\put(77,30){\makebox(0,0)[bl]{}} \put(83,30){\makebox(0,0)[bl]{}}

\put(88,29.5){\makebox(0,0)[bl]{\cite{lev06,levkud11}}}

\put(0.5,25.7){\makebox(0,0)[bl]{19.Vibration conveyor}}

\put(47,26){\makebox(0,0)[bl]{\(\star \)}}
\put(53,26){\makebox(0,0)[bl]{\(\star \)}}
\put(59,26){\makebox(0,0)[bl]{}} \put(65,26){\makebox(0,0)[bl]{}}
\put(71,26){\makebox(0,0)[bl]{}} \put(77,26){\makebox(0,0)[bl]{}}
\put(83,26){\makebox(0,0)[bl]{}}

\put(88,25.5){\makebox(0,0)[bl]{\cite{lev98}}}

\put(0.5,21.6){\makebox(0,0)[bl]{20.Geological planning}}

\put(47,22){\makebox(0,0)[bl]{\(\star \)}}
\put(53,22){\makebox(0,0)[bl]{\(\star \)}}
\put(59,22){\makebox(0,0)[bl]{}} \put(65,22){\makebox(0,0)[bl]{}}
\put(71,22){\makebox(0,0)[bl]{}} \put(77,22){\makebox(0,0)[bl]{}}
\put(83,22){\makebox(0,0)[bl]{}}

\put(88,21.5){\makebox(0,0)[bl]{\cite{lev98}}}

\put(0.5,17.5){\makebox(0,0)[bl]{21.Concrete technology}}

\put(47,18){\makebox(0,0)[bl]{\(\star \)}}
\put(53,18){\makebox(0,0)[bl]{\(\star \)}}
\put(59,18){\makebox(0,0)[bl]{}} \put(65,18){\makebox(0,0)[bl]{}}
\put(71,18){\makebox(0,0)[bl]{}} \put(77,18){\makebox(0,0)[bl]{}}
\put(83,18){\makebox(0,0)[bl]{}}

\put(88,17.5){\makebox(0,0)[bl]{\cite{lev06,levnis01}}}

\put(0.5,14){\makebox(0,0)[bl]{22.Medical treatment }}

\put(47,14){\makebox(0,0)[bl]{\(\star \)}}
\put(53,14){\makebox(0,0)[bl]{\(\star \)}}
\put(59,14){\makebox(0,0)[bl]{}} \put(65,14){\makebox(0,0)[bl]{}}
\put(71,14){\makebox(0,0)[bl]{}} \put(77,14){\makebox(0,0)[bl]{}}
\put(83,14){\makebox(0,0)[bl]{}}

\put(88,13.5){\makebox(0,0)[bl]{\cite{lev06,levsok04}}}

\put(0.5,9.5){\makebox(0,0)[bl]{23.Educational programs}}

\put(47,10){\makebox(0,0)[bl]{\(\star \)}}
\put(53,10){\makebox(0,0)[bl]{\(\star \)}}
\put(59,10){\makebox(0,0)[bl]{\(\star \)}}
\put(65,10){\makebox(0,0)[bl]{\(\star \)}}
\put(71,10){\makebox(0,0)[bl]{\(\star \)}}
\put(77,10){\makebox(0,0)[bl]{\(\star \)}}
\put(83,10){\makebox(0,0)[bl]{\(\star \)}}

\put(88,09.5){\makebox(0,0)[bl]{\cite{lev98,lev06,lev11agg,lev13}}}

\put(0.5,06){\makebox(0,0)[bl]{24.Standard in multimedia}}

\put(47,06){\makebox(0,0)[bl]{\(\star \)}}
\put(53,06){\makebox(0,0)[bl]{\(\star \)}}
\put(59,06){\makebox(0,0)[bl]{\(\star \)}}
\put(65,06){\makebox(0,0)[bl]{\(\star \)}}
\put(71,06){\makebox(0,0)[bl]{\(\star \)}}
\put(77,06){\makebox(0,0)[bl]{\(\star \)}}
\put(83,06){\makebox(0,0)[bl]{\(\star \)}}

\put(88,05.5){\makebox(0,0)[bl]{\cite{lev13mpeg,lev09had}}}

\put(0.5,01.5){\makebox(0,0)[bl]{25.Immunoassay technology}}

\put(47,02){\makebox(0,0)[bl]{\(\star \)}}
\put(53,02){\makebox(0,0)[bl]{\(\star \)}}
\put(59,02){\makebox(0,0)[bl]{}} \put(65,02){\makebox(0,0)[bl]{}}
\put(71,02){\makebox(0,0)[bl]{}} \put(77,02){\makebox(0,0)[bl]{}}
\put(83,02){\makebox(0,0)[bl]{}}

\put(88,01.5){\makebox(0,0)[bl]{\cite{lev06,levfir05}}}

\end{picture}
\end{center}

~~~


\end{document}